\documentclass[twoside,12pt]{article}
\title{Partitioning  zero-divisor graphs of finite commutative rings into global defensive alliances}
\date{}
\author{}
\usepackage[all]{xy}
\usepackage{amssymb,amsmath,latexsym}
\usepackage{theorem}
\usepackage{amsfonts}
\usepackage{ulem}
\usepackage{amscd}
\usepackage{amsxtra}
\usepackage{graphicx}
\usepackage{bm}
\usepackage{graphics} 
\usepackage{graphicx} 
\usepackage{pstricks,pst-node} 
\usepackage{subfigure}
\usepackage{tikz} 

\usepackage{tabularx} 
\usepackage{multirow} 
\usepackage{multicol} 
\usepackage{arydshln} 
\usepackage{fancybox} 
\usepackage{multicol} 
\usepackage{array} 
\usepackage{fancybox}

\usepackage{mathrsfs}
\usepackage{array,multirow,makecell}
\setcellgapes{1pt}
\makegapedcells

\newcolumntype{R}[1]{>{\raggedleft\arraybackslash }b{#1}}
\newcolumntype{L}[1]{>{\raggedright\arraybackslash }b{#1}}
\newcolumntype{C}[1]{>{\centering\arraybackslash }b{#1}}

\usetikzlibrary{arrows} 

\newcommand{\field}[1]{\mathbb{#1}}

\newcommand{\Z }{\field{Z}}

\newcommand{\F }{\field{F}}

\theoremstyle{}\newtheorem{thm}{\bf Theorem}[section]
\theoremstyle{}\newtheorem{cor}[thm]{\bf Corollary}
\theoremstyle{}
\theoremstyle{}
\theoremstyle{}
\theoremstyle{}\newtheorem{pro}[thm]{\bf Proposition}
\theoremstyle{}\newtheorem{exm}[thm]{\bf Example}
\theoremstyle{}
\theoremstyle{}
\theoremstyle{}\newcommand{\cqfd}{\hfill$\square$}

\def\pr{{\parindent0pt {\bf Proof.\ }}}
\def\cqfd
{\hspace{1cm}
	\rule{2mm}{2mm}%
	\medbreak%
	\par%
}

\def\ann{{\rm Ann}}
\def\nil{{\rm Nil}}
\begin{document}
	\thispagestyle{empty}
	
	\maketitle \vspace*{-1.5cm}
	
	
	\begin{center}
	{\large\bf   Driss Bennis$^{1, a}$, Brahim El Alaoui$^{1, b}$}
	
		\bigskip
	
		$^1$ Department of Mathematics, Faculty of Sciences,  Mohammed V University in Rabat,  Morocco.\\
		\noindent    $^a$\,driss.bennis@fsr.um5.ac.ma; driss$\_$bennis@hotmail.com;  \\ $^b$\,brahim$\_$elalaoui2@um5.ac.ma;  brahimelalaoui0019@gmail.com   \\[0.2cm]

	\end{center}
	\bigskip 
	%
	\noindent{\large\bf Abstract.}
	
For a commutative ring $R$ with identity, the zero-divisor graph of $R$, denoted $\Gamma(R)$, is the graph whose vertices are the non-zero zero divisors of $R$ with two distinct vertices $x$ and $y$ are adjacent if and only if $xy=0$. In this paper, we are interested in partitioning the vertex set of $\Gamma(R)$ into global defensive alliances for a finite commutative ring $R$. This problem has been well investigated in graph theory. Here we connected it with the ring theoretical context. We  characterize various commutative finite rings for which the zero divisor graph is partitionable into global defensive alliances. We also give several examples to illustrate the scopes and limits of our results.\\


\small{\noindent{\bf 2020 Mathematics Subject Classification :}   13M05, 05C25
		
		
		\section{Introduction}

Within this  paper, $R$  will be a commutative ring  with  $1\neq 0$,  $Z(R)$ be  its set of zero-divisors  and $U(R)$ be its set of units. Let  $x$ be an element  of  $R$, the annihilator of $x$  is defined as  $\ann(x):=\{y\in R| \ xy=0\}$. For an ideal $I$ of $R$, $\sqrt{I}$ means the radical of $I$.  An element $x$ of $R$ is called nilpotent if $x^n=0$ for some positive integers $n$. The set of all nilpotent elements is denoted   $\nil(R):=\sqrt{0}$. A ring $R$ is called reduced if $\nil(R)=\{0\}$. The ring $\Z/n\Z$ of the residues modulo an integer $n$ will be noted by $\Z_n$.  
For a subset $X$ of $R$, we denoted  $X^*=X\setminus \{0\}$.
For any real number $r$, let $\lceil r \rceil$	 (resp.,  $\lfloor r \rfloor$)
denote the ceiling of $r$, that is, the least integer greater than or equal to $r$  (resp., the  floor of $r$, that is  the greatest integer less than or equal to $r$). \\

We assume the reader has at least a basic familiarity with the   zero-divisor graph theory.  For general background   on the zero-divisor graph theory, we refer the reader to \cite{AFLL01, ADL, DM, A08, BAD,  AM07, SAAM, COY,BI}. The concept of  the zero-divisor graph of a commutative ring was first   introduced by Beck \cite{BI}, to investigate the structure of commutative rings. For a given commutative ring $R$, Beck's zero-divisor graph is a simple graph with vertex set all elements of $R$, such that two distinct vertices $x$ and $y$ are adjacent if and only if $xy=0$. Beck was mainly interested in colorings. In 1999,  Anderson and Livingston defined a simplified version $\Gamma(R)$ of Beck’s zero-divisor graph by including only nonzero zero-divisors of $R$ in the vertex set and leaving the definition of edges the same \cite{ADL}. The reason of this simplification was to better capture the essence of the zero-divisor structure of the ring. Several properties of $\Gamma(R)$ have been investigated, such as connectedness, diameter, girth, chromatic number, etc. \cite{ADL, SAAM}. In addition, the isomorphism problem for such graphs has been solved for finite reduced rings \cite{AFLL01}. Several authors have also investigated rings $R$ whose graph $\Gamma(R)$ belongs to a certain family of graphs, such as star graphs \cite{SAAM}, complete
graphs \cite{ADL}, complete $r$-partite graphs and planar graphs \cite{SHS, NOS}.\\

This paper deals with defensive alliance notions of graphs. In a graph $\Gamma$, a nonempty set of vertices $S$ is called a defensive alliance if any vertex $v$ in $S$ has at least one more neighbors in $S$ than it has in the complement of $S$   (see Section 2 for more definitions of related notions). These notions were motivated by the study of alliances between different parts in a population. Since its introduction by Kristiansen, Hedetniemi, and Hedetniemi in \cite{PSS, PKSMST}, the alliances have attracted the attention of many authors.  Recently, several authors have been  studied these notions  in the context of zero-divisor graphs of finite commutative rings (see for instance \cite{DBK, RGOJ, NA}). In this paper, we focus our attention on the problem of partitioning zero-divisor graphs of some kind of finite commutative rings into global defensive alliances. This problem has  been subject of several researches for arbitrary graphs (see for instance \cite{ Shaf, KR, ISJJ}).

This paper is organized as follows:\\
In Section 2, we recall  the global defensive alliance  of graphs as well as some notions related to it.\\
In Section 3, we study when zero-divisor graphs  of some kind of direct products of finite fields with finite
local rings is partitionable into global defensive alliances and we calculate the  global defensive alliance partition number, $\psi_g(\Gamma(R))$,  for each one of them.  We give also  complete characterizations  for partitioning  zero-divisor graphs of finite rings  with $\gamma_{a}(\Gamma(R))=1,2$.  Namely,   we prove that a zero-divisor graph, $\Gamma(R)$, of a finite ring $R$  with $\gamma_{a}(\Gamma(R))=1$ is partitionable into global defensive alliances if and only if  $R$ is isomorphic to one of the rings $\Z_2\times \Z_2$, $\Z_9$, $\Z_3[X]/(X^2)$ (see Theorem \ref{thm_partitioning1}), and  we prove that for  a zero-divisor graph of a finite ring with $\gamma_a(\Gamma(R))=2$,   $\Gamma(R)$ is partitionable into global defensive alliances if and only if $R$ is isomorphic to one of the rings $\Z_2\times \Z_4$, $\Z_2\times \Z_{2}[X]/(X^2)$, $\Z_3\times \Z_3$, $\Z_3\times \F_4$, $\Z_{25}$, $\Z_5[X]/(X^2)$ and  $\F_4\times \F_4$ (see Theorem \ref{thm_partitioning2} and Corollary \ref{cor_partition2}).


\section{Preliminaries}
In this section we deal with the alliance notion of graphs. 
We  suppose some familiarity  with some  basic concepts 
in graph theory. For the convenience of the reader, we recall this concept and some useful concepts, including the concept of the dominating set, which is a very important concept in graph theory. In fact, many interesting properties related to this term are still in the spotlight of some researchers. (see for instance \cite{HHH, TSP, TSP2}). \\

Let $\Gamma=(V,E)$ be a finite  simple graph that is a finite  graph without loop or multiple edges. Then, we will use the notation $x-y$ to mean  the edge between the   two adjacent vertices $x$ and $y$. 

Let $x\in V$ be a vertex in $\Gamma$,  the open neighborhood of $x$ is defined as $N(x):=\{y\in V| \ x-y\in E\}$,   and  the closed neighborhood of $x$ is defined by  $N[x]:=N(x)\cup \{x\}$. In general, 
for a nonempty subset $S\subseteq V$, the open neighborhood of $S$ is defined as $ N(S)=\cup_{x\in S}N(x)$ and its closed neighborhood by  $N[S]=N(S)\cup S$.\\

The degree of the vertex $x\in V$,  denoted by $deg(x)$,  is the cardinality of its open neighborhood.  Namely,    $deg(x):=|N(x)|$. In general, for every nonempty subset $S\subseteq V$ and every vertex $x\in S$,  we define the degree of $x$ over $S$ as $deg_S(x):=|S\cap N(x)|$. So, $deg_V(x)=deg(x)$.\\

A nonempty set $S\subseteq V$ is a dominating set in $\Gamma$ if for every vertex $v\in \bar{S}$, $deg_S(v)>0$. The domination  number of $\Gamma$, denoted $\gamma(\Gamma)$, is  the minimum cardinality of a dominating set  in $G$. \\

A non-empty set of vertices $S\subseteq V$  is
called a defensive alliance if for every $x\in S$, $|N[x] \cap S|\geq  |N(x) \cap \bar{S}|$,  in  other words, $deg_S(x)+1\geq deg_{\bar{S}}(x)$, where    $\bar{S}=V\setminus S$ (i.e., $\bar{S}$ is the complement of $S$ in $V$)  or equivalently $deg(v)+1\geq 2deg_{\bar{S}}(v)$.\\

A defensive alliance $S$ is called strong if for every vertex $x\in S$, $|N[x]\cap S|> |N(x)\cap \bar{S}|$,  in other words, $deg_S(x)\geq deg_{\bar{S}}(x)$.
In this case we say that every vertex in $S$ is strongly defended. A defensive alliance $S$ is global if it forms a dominating set.\\

The global defensive alliance partition number of $\Gamma$, denoted $\psi_g(\Gamma)$, is defined to be the maximum number of sets in a partition of $V$ such that each set is a global defensive alliance. A graph $\Gamma$ is partitionable into global defensive alliances if $\psi_g(\Gamma)\geq 2$.

		
\section{Partitioning  into  global defensive alliances}
In this section we study when  zero-divisor graphs of some kind of commutative rings can be partitinable into global defensive alliances and calculate their global defensive alliance partition numbers.
First we start by   giving  a lower bound of  the cardinality of the set of zero divisors $Z(R)$ in term of the global defensive alliance partition number of $\Gamma(R)$.

\begin{pro}\label{prop_|Z(R)|}
Let $R$ be a ring. Then,  $$ |Z(R)|\geq \psi_g(\Gamma(R))^2-\psi_g(\Gamma(R))+1.$$ 
\end{pro}
	 
\pr
If $\Gamma(R)$ is not partitionable into global defensive alliances, then $\psi_g(\Gamma(R))=1$ and so $|Z(R)|\geq \psi_g(\Gamma(R))^2-\psi_g(\Gamma(R))+1=1$ it is true. So,  we  assume that $\Gamma(R)$ is partitionable into global defensive alliances. Let $r=\psi_g(\Gamma(R))$ and $\{S_1,\ldots, S_r\}$ be a partition of $\Gamma(R)$ into global defensive alliances.  Since, for every $i\in [\![1;r]\!]$, $S_i$ is a dominating set,  $deg_{\bar{S_i}}(x)\geq r-1$ for every $x\in S_i$ and so  $|S_i|-1\geq \deg_{S_i}(x)\geq deg_{\bar{S_i}}(x) -1$ (since $S_i$ is a defensive alliance). Thus, $|S_i|-1\geq r-2$ and hence $|Z(R)|-1=\sum_{i=1}^{r}|S_i|\geq r^2-r$. Then, $|Z(R)|\geq r^2-r+1$. On the other hand, if $\Gamma(R)$ is not partitionable into global defensive alliances,  that is $\psi_g(\Gamma(R))=1$, we have  $|Z(R)|\geq 1^2-1+1=1$. Hence, $ |Z(R)|\geq \psi_g(\Gamma(R))^2-\psi_g(\Gamma(R))+1$.
\cqfd

The following example shows that this bound is sharp.

\begin{exm}
			Let $R=\Z_9$. Then, the  zero-divisor graph is just an edge joining $\bar{3}$ and $\bar{6}$. So, $\{S_1,S_2\}$ where   $S_1=\{\bar{3}\}$ and $S_2=\{\bar{6} \}$, is a partition of $\Gamma(R)$ into two  global defensive alliances. Then, $\psi_g(\Gamma(R))=2$ and so $|Z(R)|= \psi_g(\Gamma(R))^2-\psi_g(\Gamma(R))+1=2^2-2+1=3$.
\end{exm}

For a finite local ring $(R,M)$ which is not a field, $M=Z(R)=\ann(x)$ for some $x\in Z(R)^*$ and  $|R|=p^{nr}$ and $|M|=p^{(n-1)r}$ for some prime number $p$ and positive integers $n$ and $r$. However, we know that for a finite local ring $(R,M)$, $\Gamma(R)$ is complete if and only if $Z(R)=M$ with $M^2=0$,   \cite[Theorem 2.8]{ADL}.  So, we have the following result for this simple case.

\begin{pro}
	Let $(R,M)$ be a finite local  ring such that its maximal ideal $M$ is nilpotent of index $2$. Then,  
	      \begin{enumerate}
	      	    \item  if $|M|$ is odd,  $\Gamma(R)$ is  partitionable into global defensive alliances with  $\psi_g(\Gamma(R))=2$.
	      	    \item if $|M|$ is even,  $\Gamma(R)$ is not partitionable into global defensive alliances.
	      \end{enumerate}
\end{pro}

\pr 
  (1)- Let    $S_1$ and $S_2$ be two distinct subsets of $M^*$ such that $|S_1|=|S_2|= \frac{|M|-1}{2}$. Then, $\{S_1,S_2\}$ is a partition of $\Gamma(R)$ into global defensive alliances and so $\psi_g(\Gamma(R))\geq 2$. Since, $\gamma_{a}(\Gamma(R))= \frac{|M|-1}{2}$ and  $\gamma_{a}(\Gamma(R))\times\psi_g(\Gamma(R))\leq |M|-1$,  $\psi_g(\Gamma(R))\leq 2$. Hence, $\psi_g(\Gamma(R))=2$.\\
 (2)- Suppose that $\Gamma(R)$ is partitionable into global defensive alliances, then $\psi_g(\Gamma(R))\geq 2$ and so, by \cite[Proposition 3.6]{DBK},  $\left\lceil\frac{|M|-1}{2}\right\rceil\times 2\leq \gamma_a(\Gamma(R))\times\psi_g(\Gamma(R))\leq |M|-1$. Thus, $|M|\leq |M|-1$, a contradiction. 
\cqfd

\begin{cor}
	Let $p$ be a prime number. Then,  we have two cases:
	\begin{enumerate}
		\item if $p=2$, then $\Gamma(\Z_{p^2})$ has only one vertex  and so it is not partitionable into global defensive alliances.
		\item  if $p\neq 2$, then  $\Gamma(\Z_{p^2})$ is partitionable into global defensive alliances and $\psi_g(\Gamma(\Z_{p^2}))=2$.
	\end{enumerate}
\end{cor}

In the following theorem we study when $\Gamma(\Z_{p^n})$ is partitionable into global defensive alliances for a prime number $p$ and a positive integer $n\geq 3$.

\begin{thm}
    Let $p$ be a prime number and $n\geq 3$ be a positive integer. Then,
    \begin{enumerate}
    	\item If $p=2$, then $\Gamma(\Z_{p^n})$ is not partitionable into  global defensive alliances.
    	\item If $p\geq 3$, then  $\Gamma(\Z_{p^n})$ is partitionable into  global defensive alliances and $\psi_{g}(\Gamma(\Z_{p^n}))=2$. 
    \end{enumerate}
\end{thm}
\pr 
We have $Z:=Z(\Z_{p^n})=\{\overline{mp}|\  0\leq m < p^{n-1}\}$  and $|Z(\Z_{p^n})|=p^{n-1}$. Then,\\
(1)- Suppose that $\Gamma(\Z_{p^n})$ is partitionable into global defensive alliances. Then, $\psi_{g}(\Gamma(\Z_{p^n}))\geq 2$. Since $\gamma_{a}(\Gamma(\Z_{p^n}))=2^{n-2}$, by \cite[Theorem 2.9 ]{NA}, then  $2^{n-2}\times 2\leq \gamma_a(\Gamma(\Z_{p^n}))\psi_{g}(\Gamma(\Z_{p^n}))\leq 2^{n-1}-1 $, a contradiction.\\
(2)-For each $1\leq k \leq n-1$, set $A_k = \{\overline{ap^k}\in Z |\ p \text{ does not divide } a\}$. The sets $A_k$ are disjoints, $|A_k|=p^{n-k}-p^{n-k-1}$ which is an even number, and $Z^*=\cup_{k=1}^{n-1} A_k$.
Let $S_1=\cup_{k=1}^{n-1} A_k'$ and $S_2=\cup_{k=1}^{n-1}  A_k''$ such that $A_k'$ is a half of elements of   $A_k$ and $A_k''$ the other half. By the proof of \cite[Theorem 2.9]{NA}, $S_1$ and $S_2$ are two global defensive alliances and so $\{S_1,S_2\}$ is a partition of $\Gamma(\Z_{p^n})$ into global defensive alliances. Then, $\psi_{g}(\Gamma(\Z_{p^n}))\geq 2$. On the other hand $\gamma_a(\Gamma(\Z_{p^n}))=\left\lceil\frac{p^{n-1}-1}{2}\right\rceil$,  by \cite[Theorem 2.9 ]{NA},  and since  $\gamma_a(\Gamma(\Z_{p^n}))\psi_{g}(\Gamma(\Z_{p^n}))\leq p^{n-1}-1 $,   $\psi_{g}(\Gamma(\Z_{p^n}))\leq 2$. Hence, $\psi_{g}(\Gamma(\Z_{p^n}))= 2$
\cqfd

The following result  characterizes when a zero-divisor graph of $\Z_2 \times F$, for a finite field $F$, is partitionable into global defensive alliances.  	
	
\begin{thm}
	Let $F$ be a finite field. Then,  $\Gamma(\Z_2\times F)$ is partitionable into global defensive alliances if    and only if  $F\cong \Z_2$.
\end{thm}		
\pr
$\Leftarrow)$ Let  $S_1=\{(1,0)\}$ and $S_2=\{(0,1)\}$, then $S_1$ and $S_2$ are both global defensive alliances and so  $\{S_1,S_2\}$ is the only partition into global defensive alliances of $\Gamma(\Z_2\times \Z_2)$. Then, $\psi_g(\Gamma(\Z_2\times \Z_2))=2$.\\
$\Rightarrow)$ Assume that $|F|\geq 3$ and suppose  that $\Gamma(\Z_2\times F)$ is partitionable into global defensive alliances. Let $S_1$ be  a global defensive alliance in a partition of $\Gamma(\Z_2\times F)$ such  $(1,0)\notin S_1$, then $\{0\}\times F^*\subset S_1$ (since, $S_1$ is a dominating set) and so $S_2=\{(1,0)\}$ is the other global defensive alliance such that $\{S_1, S_2\}$ is a partition of $\Gamma(\Z_2\times F)$ into global defensive alliances. Then,  $deg_{S_2}(1,0)+1\geq deg_{\bar{S_2}}(1,0)$  and so $1\geq |F|-1$, a contradiction.
\cqfd

\begin{cor}
   Let $p$ be a prime number. Then, 
	\begin{enumerate}
	\item If $p=2$, then $\Gamma(\Z_2\times \Z_p)$ is partitionable into global defensive alliances   and $\psi_g(\Gamma(\Z_2\times \Z_p))=2$.
	\item  If $p\neq 2$, then $\Gamma(\Z_2\times \Z_p)$ is not partitionable into global defensive alliances.
\end{enumerate}    
\end{cor}
 The  following theorem  characterizes when $\Gamma(\Z_2\times R)$, for a finite local ring $R$ which is not a field,  is partitionable into global defensive alliances.  
\begin{thm}
	Let $R$ be a finite local ring which is not a field. Then,
	$\Gamma(\Z_2\times R)$ is partitionable into global defensive alliances  if and only if $R\cong \Z_4$ or $R\cong \Z_2[X]/(X^2)$.	\\
	
	Moreover, if $\Gamma(\Z_2\times R)$ is partitionable into global defensive alliances, then  $\psi_g(\Gamma(\Z_2\times R))=2$.
\end{thm}
\pr
$\Leftarrow)$ Assume that $R\cong \Z_4$. The zero-divisor graph of $\Z_2\times \Z_4$ is illustrated in Figure \ref{fig2}.

\begin{figure}[ht]
	\centering
	\includegraphics[scale=0.5]{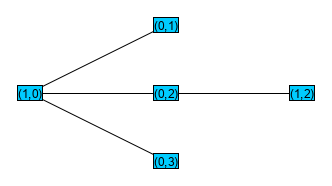}
	\caption{$\Gamma(\Z_2\times \Z_4)$}
	\label{fig2}
\end{figure}
 Set $S_1=\{(1,0),(0,2)\}$ and $S_2=\{(0,1),(0,3),(1,2)\}$. So, $S_1$ and $S_2$ are two global defensive alliances since $S_1$ and $S_2$ are dominating sets and

 	\begin{equation*}
 \left\{
 \begin{array}{ll}
 & deg_{S_1}((\bar{1},\bar{0}))+1=2\geq deg_{\bar{S_1}}((\bar{1},\bar{0}))=2,\\&
 deg_{S_1}((\bar{0},\bar{2}))+1=2\geq deg_{\bar{S_1}}((\bar{0},\bar{2}))=1,\\&
 deg_{S_2}((\bar{0},\bar{1}))+1=1\geq deg_{\bar{S_2}}((\bar{0},\bar{1}))=1,\\&
 deg_{S_2}((\bar{0},\bar{3}))+1=1\geq deg_{\bar{S_2}}((\bar{0},\bar{3}))=1,\\&
 deg_{S_2}((\bar{1},\bar{2}))+1=1\geq deg_{\bar{S_2}}((\bar{1},\bar{2}))=1.
 \end{array}
 \right.
 \end{equation*}
 On the other hand $Z(\Z_2\times \Z_4)^*=S_1\cup S_2$ and $S_1\cap S_2\neq \emptyset$ and so $\{S_1,S_2\}$ form a partition of $\Gamma(\Z_2 \times R)$ into global defensive alliances. Then,  $\psi_g(\Gamma(\Z_2 \times R))\geq 2$. Thus, by  \cite[Theorem 2.1 ]{ISJJ},  $\psi_g(\Gamma(\Z_2 \times R))\leq \left \lfloor \frac{\delta+1+2}{2}\right\rfloor=2$  and so  $\psi_g(\Gamma(\Z_2 \times R))=2$. Similarly, when $R\cong \Z_2[X]/(X^2)$, we take the partition $\{S_1,S_2\}$  with $S_1=\{(\bar{1},\bar{0}),(\bar{0},\bar{X})\}$  and $S_2=\{(\bar{0},\bar{1}),(\bar{0},\overline{1+X}),(\bar{1},\bar{X})\}$.

$\Rightarrow)$ Assume that $R\ncong \Z_4$ and $R\ncong \Z_2[X]/(X^2)$. Suppose that $\Gamma(\Z_2\times R)$ is  partitionable into   two global defensive   alliances $S_1$ and $S_2$. So, assume that $(1,0)\in S_1$. Then, we have two cases:\\
Case 1: $\{0\}\times U(R)\subset S_2$ and so $ deg_{S_1}((\bar{1},0))+1\geq deg_{\bar{S_1}}((\bar{1},0))$. Then,  $|Z(R)|-1+1\geq deg_{S_1}((\bar{1},0))+1\geq deg_{\bar{S_1}}((\bar{1},0))\geq |U(R)|$, a contradiction.\\
Case 2:there exists $u\in U(R)$ such that  $(0,u)\in S_1$. Then, $S_2$ is not a dominating set since there is no vertex adjacent to $(0,u)$ other than $(1,0)$, a contradiction.\\
Hence, $\Gamma(\Z_2\times R)$ is not   partitionable into   global defensive alliances.
\cqfd

\begin{cor}
  Let $p$ be a prime number and $n\geq 2$ be a positive integer. Then, $\Gamma(\Z_2\times \Z_{p^n})$ is partitionable into global defensive alliances if and only if $p=n=2$. 
\end{cor}

In the following result, we study when zero-divisor graphs of a direct product of two finite fields is partitionable into  global defensive alliances.

\begin{thm}
	Let $F$ and $K$ be two finite fields such that $|F|,|K|\geq 3$. Then,  $\Gamma(F \times K)$ is partitionable into global defensive alliances and 
		
 \begin{equation*}
		\psi_g(\Gamma(F\times K))=  \left\{
	\begin{array}{ll}
	& 	3 \text{\hspace{2cm} if } |F|=|K|=4,\\&
	2 \text{\hspace{2cm} otherwise. }
	\end{array}
	\right.
\end{equation*}
\end{thm}

\pr We have the following cases:\\
\textbf{Case} $|F|=|K|=4$: The zero-divisor graph of $F\times K$ is illustrated  in Figure \ref{fig3}.
 \begin{figure}[ht]
	\centering
	\includegraphics[scale=0.5]{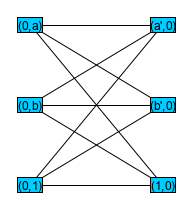}
	\caption{$\Gamma(F\times K)$}
	\label{fig3}
\end{figure}

 So, $\{\{(0,a), (a',0)\}, \{(0,b),(b',0)\}, \{(0,1),(1,0)\} \}$ is a partition of  $\Gamma(F\times K)$  into  global defensive alliances, then $\psi_g(\Gamma(F\times K))\geq 3$. On the other hand,  $\gamma_a(\Gamma(F\times K))\psi_g(\Gamma(F\times K))\leq 6$. Thus, $\psi_g(\Gamma(F\times K)) =3$.\\ 
\textbf{Case} $|F|\neq 4$ or $|K|\neq 4$: 
Let $S_1=F_1\times \{0\}\cup \{0\}\times K_1$ and $S_2=F_2\times \{0\}\cup \{0\}\times K_2$ such that  $F_1, F_2\subset F^*$ and $K_1,K_2\subset K$ with 

 	\begin{equation*}
\left\{
\begin{array}{ll}
& F_1\cap F_2=K_1\cap K_2=\emptyset,\\&
|F_1|=\left \lfloor \frac{|F|-1}{2}\right \rfloor,\\&
|F_2|= |F|-1- \left \lfloor \frac{|F|-1}{2}\right \rfloor,\\&
|K_1|= \left \lfloor \frac{|K|-1}{2}\right \rfloor,\\&
|K_2|= |K|-1- \left \lfloor \frac{|K|-1}{2}\right \rfloor.
\end{array}
\right.
\end{equation*}

So, $S_1$ and $S_2$ are two global defensive alliances and $Z(F\times K)^*=S_1\cup S_2$ and $S_1\cap S_2\neq \emptyset$. Thus,   $\Gamma(F \times K)$ is partitionable into global defensive alliances and so $\psi_g(\Gamma(F\times K))\geq 2$.\\
Now, suppose that $\psi_g(\Gamma(F\times K))\geq 3$, then  $3\times \gamma_a(\Gamma(F\times K))\leq \gamma_a(\Gamma(F\times K))\psi_g(\Gamma(F\times K))\leq |F|+|K|-2$ and since $\gamma_a(\Gamma(F\times K))=\left \lfloor \frac{|F|-1}{2}\right \rfloor+\left \lfloor \frac{|K|-1}{2}\right \rfloor$, by \cite[Proposition 2.3]{NA}, then  $3\times (\left \lfloor \frac{|F|-1}{2}\right \rfloor+\left \lfloor \frac{|K|-1}{2}\right \rfloor) \leq |F|+|K|-2$. So, we have four sub-cases to discuss:\\
\textbf{sub-case 1}; $2\mid_{|F|}$ and $2\mid_{|K|}$: Then, $|F|+|K|\leq 8$,  a contradiction since one of the $|F|$ and $|K|$ is different from $4$.\\
\textbf{sub-case 2}; $2\mid_{|F|}$ and $2\nmid_{|K|}$: Then, $|F|+|K|\leq 5$,  a contradiction.\\
\textbf{sub-case 3}; $2\nmid_{|F|}$ and $2\mid_{|K|}$: Similar to sub-case 2.\\
\textbf{sub-case 4}; $2\nmid_{|F|}$ and $2\nmid_{|K|}$: Then, $|F|+|K|\leq 2$, a contradiction.\\
Hence, $\psi_g(\Gamma(F\times K))=2$.
\cqfd
\begin{cor}
	Let $p, q\geq 3$ be two prime numbers. Then, $\Gamma(\Z_p\times \Z_q)$ is partitionable into global defensive alliances and  
$\psi_g(\Gamma(\Z_p\times \Z_q))=2$.
\end{cor}

\begin{thm}\label{thm_F times R}
	Let $R$ be a finite local  ring such that its maximal ideal  is nilpotent of index $2$ and  $F$ be a finite field with $|F|\geq 3$.  Then,
\begin{enumerate}
	\item If $|Z(R)|$ is odd, then $\Gamma(F\times R)$ is partitionable into global defensive alliances and $\psi_g(\Gamma(F\times R))=2$.
	\item If $|Z(R)|$ is even, then  $\Gamma(F\times R))$ is not partitionable into global defensive alliances.  
\end{enumerate}	
\end{thm}

\pr (1)-  Assume that  $|Z(R)|$ is odd and set $S_1=A_1\times\{0\}\cup  \{0\} \times A_2\cup \{0\}\times A_3\cup A_1\times A_3$ 
and $S_2=B_1\times\{0\}\cup \{0\}\times B_2\cup \{0\}\times B_3\cup B_1\times B_3$ such that

\begin{equation*}
\left\{
\begin{array}{ll}
& A_1, B_1\subset F^* \text{ with } |A_1|=\left \lceil \frac{|F|-1}{2}\right \rceil, \ |B_1|=|F|-1-|A_1| \text{ and  }  A_1\cap B_1=\emptyset,\\&
A_2, B_2\subset U(R) \text{ with } |A_2|=\left \lceil \frac{|U(R)|}{2}\right \rceil, \ |B_2|=|U(R)|-|A_2| \text{ and  }  A_2\cap B_2=\emptyset,\\&
A_3, B_3\subset Z(R)^* \text{ with } |A_3|=\left \lceil \frac{|Z(R)|-1}{2}\right \rceil, \ |B_3|=|Z(R)|-1-|A_3| \text{ and  }  A_3\cap B_3=\emptyset.
\end{array}
\right.
\end{equation*} 
It is clear that $S_1$ is a dominating set. So, let prove that $S_1$ is a defensive alliance. Let $(x,0)\in A_1\times\{0\}$, then  $deg_{S_1}((x,0))+1= |A_2|+|A_3|+1=\left \lceil \frac{|U(R)|}{2}\right \rceil + \left \lceil \frac{|Z(R)|-1}{2}\right \rceil +1$ and $deg_{\bar{S_1}}((x,0))= |B_2|+|B_3|=|R|-1-(\left \lceil \frac{|U(R)|}{2}\right \rceil + \left \lceil \frac{|Z(R)|-1}{2}\right \rceil)$ and so $deg_{S_1}((x,0))+1\geq deg_{\bar{S_1}}((x,0))$. Let $(0,y)\in \{0\}\times A_2$, then $deg_{S_1}((0,y))+1=|A_1|=\left \lceil \frac{|F|-1}{2}\right \rceil +1$ and $deg_{\bar{S_1}}((0,y))=|B_1|=|F|-1-\left \lceil \frac{|F|-1}{2}\right \rceil\leq \left \lceil \frac{|F|-1}{2}\right \rceil+1$ and so $deg_{S_1}((0,y))+1\geq deg_{\bar{S_1}}((0,y))$. Let $(0,y)\in \{0\}\times A_3$, then $deg_{S_1}((0,y))+1=|A_1|+|A_3|-1+|A_1||A_3|+1=|A_1|+|A_3|+|A_1||A_3|$ and  $deg_{\bar{S_1}}((0,y))=|B_1|+|B_3|+|B_1||B_3|=|F||Z(R)|-|F||A_3|-|A_1||Z(R)|+|A_1||A_3|-1$. Then, if $2$ divide $|F|$, then $deg_{S_1}((0,y))=\frac{|F|}{2}+\frac{|Z(R)|-1}{2}+|A_1||A_3|$ and  $deg_{\bar{S_1}}((0,y))=\frac{|F|}{2}+|A_1||A_3|-1$ (since $|Z(R)|$ is odd) and so   $deg_{S_1}((0,y))+1\geq deg_{\bar{S_1}}((0,y))$, otherwise $deg_{S_1}((0,y))=\frac{|F|-1}{2}+\frac{|Z(R)|-1}{2}+|A_1||A_3|$ and   $deg_{\bar{S_1}}((0,y))=\frac{|F|}{2}+\frac{|Z(R)|}{2}+|A_1||A_3|-1$ (since $|Z(R)|$ is odd) and so   $deg_{S_1}((0,y))+1\geq deg_{\bar{S_1}}((0,y))$. Finally, let $(x,y)\in A_1\times A_3$, then $deg_{S_1}((x,y))+1=|A_3|+1= \left \lceil \frac{|Z(R)|-1}{2}\right \rceil +1=\frac{|Z(R)|+1}{2}$ (since $|Z(R)|$ is odd) and  $deg_{\bar{S_1}}((x,y))=|B_3|=|Z(R)|-1-|A_3|=\frac{|Z(R)|-1}{2}$ (since $|Z(R)|$ is odd), then $deg_{S_1}((x,y))+1\geq deg_{\bar{S_1}}((x,y))$.  Then, $S_1$  is a global defensive alliance. Similarly, we prove that $S_2$ is a global defensive alliance. Since $S_1\cap S_2=\emptyset$ and $S_1\cup S_2=Z(F\times R)^*$, $\{S_1,S_2\}$ is a partition of $\Gamma(F\times R)$ into global defensive alliances. Thus, $\psi_g(\Gamma(F\times R))\geq 2$. Now, suppose that  $\psi_g(\Gamma(F\times R))\geq 3$. Then, $3\times \gamma_a(\Gamma(F\times R))\leq |Z(F\times R)|-1$ and so by \cite[Theorem 2.5 ]{NA} and since $|Z(R)|$ is odd,  $1\leq \frac{|R|}{2}+\frac{(|F|-2)|Z(R)-1|}{2}+|Z(R)|\leq 0$, a contradiction. Hence, $\psi_g(\Gamma(F\times R))=2$.\\
(2)- Suppose  that $\Gamma(F\times R)$ is partitionable into two global defensive alliances, $S_1$ and $S_2$. There are two cases to discuss:\\
 \textbf{Case} $|Z(R)|=2$ (assume that $R\cong \Z_4$):   Then, let  assume that $(0,\bar{2})\in S_1$, then $F^*\times \{\bar{2}\}\subset S_2$ (since $S_2$ is a dominating set). Since $deg_{S_2}((0,\bar{2}))+1\geq deg_{\bar{S_2}}((0,\bar{2}))$, $S_1$ contains at least $|F^*|-1$ elements from $F^*\times \{\bar{0}\}$. If $F^*\times \{\bar{0}\}\subset S_1$, then either $(0,\bar{1})$ and $(0,\bar{3})$ are both in $S_1$ or one of them is in $S_2$, if $\{(0,\bar{1}), (0,\bar{3})\}\subset S_1$, then $S_2$ is not a dominating set, a contradiction, if one of them (i.e., $(0,\bar{1})$ or $(0,\bar{3})$ ), say $(0,\bar{1})$, is in $S_2$, then $deg_{S_2}((0,\bar{1}))+1\geq deg_{\bar{S_2}}((0,\bar{1}))$ and so $|F|\leq 2$, a contradiction. Hence, there exists $(x,\bar{0})\in S_2$ for some $x\in F^*$ and so $deg_{S_1}((0,\bar{2}))+1\geq deg_{\bar{S_1}}((0,\bar{2}))$, which implies $|F^*|-1+1\geq |F^*|+1$, a contradiction. Thus,  $\Gamma(F\times R)$ is not  partitionable into  global defensive alliances.\\
 \textbf{Case} $|Z(R)|>2$: Suppose that $\{0\}\times Z(R)^*\subset S_1$, then $F^*\times Z(R)^*\subset S_2$ (Since $S_2$ is dominating set). Thus, for every $(x,r)\in F^*\times Z(R)^*\cap S_2$, $deg_{S_2}((x,r))+1\geq deg_{\bar{S_2}}((x,r))$ and so $|Z(R)|\leq 2$, a contradiction. Then, $\{0\}\times Z(R)^*\cap S_2\neq \emptyset $ (similarly, $\{0\}\times Z(R)^*\cap S_1 \neq \emptyset$). Analogously, if $F^*\times Z(R)^* \subset S_2$, then for every $(0,r)\in S_1\cap \{0\}\times Z(R)^*$, $|F|-1+|Z(R)|-3+1\geq deg_{S_1}((0,r))+1\geq deg_{\bar{S_1}}((0,r))\geq |F^*||Z(R)^*|+1$ and so $|F|+|Z(R)|\geq \frac{|F||Z(R)|+1}{2}+2$, a contradiction. Thus, $F^*\times Z(R)^*\cap S_2 \neq \emptyset$ and $F^*\times Z(R)^*\cap S_1\neq \emptyset$. Since, for every $(x,r)\in F^*\times Z(R)^*\cap S_2$ and $(x',r')\in F^*\times Z(R)^*\cap S_1$, $deg_{S_2}((x,r))+1\geq deg_{\bar{S_2}}((x,r))$ and  $deg_{S_1}((x',r'))+1\geq deg_{\bar{S_1}}((x',r'))$, then we can assume that $|S_1\cap \{0\}\times Z(R)^*|=\frac{|Z(R)|}{2}$ and $|S_2\cap \{0\}\times Z(R)^*|=\frac{|Z(R)|}{2}-1$. Now, suppose that $\{0\}\times U(R)\subset S_1$. Then, there exists a $(x,0)\in S_2$ and so $deg_{S_2}((x,0))+1\geq deg_{\bar{S_2}}((x,0))$ which implies that $\frac{|Z(R)|}{2}-1+1\geq \frac{|Z(R)|}{2} +|U(R)|$, a contradiction. So, $\{0\}\times U(R)\cap S_1\neq \emptyset $ and $\{0\}\times U(R)\cap S_2\neq \emptyset$. Then, there exist $(0,u)\in \{0\}\times U(R)\cap S_1$ and $(0,u')\in \{0\}\times U(R)\cap S_2$ such that $deg_{S_1}((0,u))+1\geq deg_{\bar{S_1}}((0,u))$ and $deg_{S_2}((0,u'))+1\geq deg_{\bar{S_2}}((0,u'))$. Then, $|F^*\times \{0\}\cap S_1|=\left \lceil \frac{|F|-1}{2}\right\rceil$ and $|F^*\times \{0\}\cap S_2|=|F|-1-\left \lceil\frac{|F|-1}{2}\right\rceil$ (or $|F^*\times \{0\}\cap S_2|=\left \lceil \frac{|F|-1}{2}\right\rceil$ and $|F^*\times \{0\}\cap S_1|=|F|-1-\left \lceil\frac{|F|-1}{2}\right\rceil$ )    and so for every $(0,r)\in S_1\cap \{0\}\times Z(R)^*$ and $(0,r')\in S_2\cap \{0\}\times Z(R)^*$  one of the following inequalities does not hold: $deg_{S_1}((0,r))+1\geq deg_{\bar{S_1}}((0,r))$ and $deg_{S_2}((0,r'))+1\geq deg_{\bar{S_2}}((0,r'))$, a contradiction. Hence,  $\Gamma(F\times R))$ is not partitionable into global defensive alliances. 
 \cqfd

\begin{cor}
	Let $p$ and $q$ be two prime numbers such that $p\neq 2$. Then, $\Gamma(\Z_p\times \Z_{q^2})$ is partitionable into global defensive alliances if and only if $q\neq 2$. Namely,
	
	 \begin{equation*}
	\psi_g(\Gamma(\Z_p\times \Z_{q^2}))=  \left\{
	\begin{array}{ll}
	& 	2 \text{\hspace{2cm} if } q\neq 2,\\&
	1 \text{\hspace{2cm} otherwise. }
	\end{array}
	\right.
	\end{equation*}
\end{cor}

To give an  example for the second assertion in Theorem \ref{thm_F times R}, in the case $|Z(R)|>2$,  we will use the idealization. Recall  that  the idealization of an $R-$module $M$ called also the trivial extension of $R$ by $M$, denoted by $R(+)M$, is the commutative ring $R\times M$ with the following addition and multiplication: $(a,n)+(b,m)=(a+b,n+m)$ and $(a,n)(b,m)=(ab,am+bn)$  for every   $(a,n),(b,m)\in  R(+)M$.

\begin{exm}
	Let  $n\geq 2$ be a positive integer, and $p$ and $q$ be two prime numbers. Then, $\Z_q(+)(\Z_{q})^n$ is a finite local  ring of maximal ideal  $0(+)(\Z_{q})^n.$  We have $(0(+)(\Z_{q})^n)^2=0$ and  so 
	$\Gamma(\Z_p \times (\Z_q(+)(\Z_{q})^n))$ is not partitionable into global defensive alliances if and only if $n\times q$ is even. Namely,
	
		 \begin{equation*}
	\psi_{g}(\Gamma(\Z_p \times (\Z_q(+)(\Z_{q})^n)))=\left\{
	\begin{array}{ll}
	& 	2 \text{\hspace{2cm} if } q\neq 2 \ \  \text{and}  \ n \text{ is odd } ,\\&
	1 \text{\hspace{2cm} otherwise. }
	\end{array}
	\right.
	\end{equation*}
\end{exm}

\begin{thm}
	Let $F$ be a finite field. Then, 
	\begin{enumerate}
		\item If $|F|\leq 4$, Then, $\Gamma(\Z_2\times \Z_2 \times F)$ is partitionable  into global defensive alliances with  $\psi_g(\Gamma(\Z_2\times \Z_2\times F))=2$.
		\item  If $|F|>4$, then $\Gamma(\Z_2 \times \Z_2 \times F )$ is not partitionable into global defensive alliances. 
	\end{enumerate} 
\end{thm}
\pr (1)-If $|F|=2$ that is $F\cong \Z_2$, then we take the partition $\{S_1,S_2\}$ such that $S_1=\{(1,0,0),(0,0,1),(0,1,0)\}$ and $S_2=\{(1,1,0),(0,1,1), (1,0,1)\}$. It is easy to see that $S_1$ and $S_2$ are both global defensive alliances. Since $Z(\Z_2\times\Z_2\times F)=\{(0,0,0)\}\cup S_1\cup S_2$, $\Gamma(\Z_2\times \Z_2 \times F)$ is partitionable into global defensive alliances and $\psi_g(\Gamma(\Z_2\times \Z_2\times \Z_2))=2$. Otherwise, we take the partition $\{S_1,S_2\}$ such that $S_1=\{(1,0,0), (0,1,0)\}\cup \{(0,0,x)|\ x\in F-\{0,1\}\}$ and $S_2=\{(0,0,1), (1,1,0)\}\cup \{(0,1,x)|\ x\in F^*\} \cup \{(1,0,x)| \ x\in F^* \}$. Hence, $\Gamma(\Z_2\times \Z_2 \times F)$ is partitionable into global defensive alliances and $\psi_g(\Gamma(\Z_2\times \Z_2\times F))=2$. \\
(2)-$\Gamma(\Z_2\times \Z_2\times F)$ is a simple graph of minimal degree $\delta= 1$, then by \cite[Theorem 2.1]{ISJJ}, $\psi_g(\Gamma(\Z_2\times \Z_2 \times F))\leq 2$. Suppose that $\{S_1,S_2\}$ is a partition of  $\Gamma(\Z_2\times \Z_2 \times F)$ into global defensive alliances.
 Assume that $(0,1,0)\in S_1$. 
 Then, we have the following two cases:\\
\textbf{Case 1} $(1,0,0)\in S_1$: If  $\{0\}\times\{0\}\times F^*\subset S_1$, then $(1,1,0)\in S_2$ (since $S_2$ is a dominating set) and so $deg_{S_2}((1,1,0))+1\geq deg_{\bar{S_2}}((1,1,0))$, then $2\geq |F|$, a contradiction. Otherwise, there exists a vertex $(0,0,u)\in S_2$, then $(1,1,0)\in S_2$, otherwise $deg_{S_2}((0,0,u))+1\geq deg_{\bar{S_2}}((0,0,u))$ implies $1\geq 3$, a contradiction. On the other hand $\{0\}\times \{1\}\times F^*\subset S_2$ (since $S_2$ is a dominating set). Thus, if there exists $u\neq v\in F^*$ such that $(0,0,v)\in S_2$, then $deg_{S_1}((1,0,0))+1\geq deg_{\bar{S_1}}((1,0,0))$ and so $-1\geq 0$, a contradiction, otherwise $deg_{S_2}((1,1,0))+1\geq deg_{\bar{S_2}}((1,1,0))$ implies that $4\geq |F|$, a contradiction.\\       
\textbf{Case 2} $(1,0,0)\in S_2$: If $\{0\}\times \{0\} \times F^*\subset S_2$, then $(1,1,0)\in S_1$ (since $S_1$ is a dominating set) and so $deg_{S_1}((1,1,0))+1\geq deg_{\bar{S_1}}((1,1,0))$, thus $2\geq |F|$, a contradiction. Then, there exists $u\in F^*$ such that  $(0,0,u)\in S_1$ and so $deg_{S_2}((0,1,0))+1\geq deg_{\bar{S_2}}((0,1,0))$  which implies that $0\geq 2$, a contradiction. 

Hence, $\Gamma(\Z_2 \times \Z_2 \times F )$ is not partitionable into global defensive alliances.
\cqfd

\begin{cor}
	Let $p$ be a prime number. Then, $\Gamma(\Z_2\times \Z_2\times \Z_p)$ is partitionable  into global defensive alliances if and only if  $p\in \{2,3\}$. Namely,  
	
	 \begin{equation*}
	\psi_g(\Gamma(\Z_2 \times \Z_2 \times \Z_p))=  \left\{
	\begin{array}{ll}
	& 	2 \text{\hspace{2cm} if }\  p\in \{2,3\},\\&
	1 \text{\hspace{2cm} otherwise.}
	\end{array}
	\right.
	\end{equation*}
\end{cor}

\begin{thm}
	Let $F$ and $K$ be two finite fields such that $|F|,|K|\geq 3$. Then, $\Gamma(\Z_2\times F\times K)$ is not partitionable into global defensive alliances.
\end{thm}
\pr 
$\Gamma(\Z_2 \times F \times K)$ is a simple graph of minimal degree $\delta=1$. Then, by \cite[Theorem 2.1]{ISJJ}, $\psi_g(\Gamma(\Z_2\times F\times K))\leq 2$. Suppose that $\psi_g(\Gamma(\Z_2\times F\times K))= 2$ that is  $\Gamma(\Z_2 \times F \times K)$ is  partitionable  into two global defensive alliance, say $S_1$ and $S_2$.  Assume that $(1,0,0)\in S_1$, then $\{0\}\times F^*\times K^*\subset S_2$. Thus, $deg_{S_1}((1,0,0))+1\geq deg_{\bar{S_1}}((1,0,0))$ implies that $|F^*|+|K^*|+1\geq |F^*||K^*|$ (since $deg_{S_1}((1,0,0))\leq |F^*|+|K^*|$ and $deg_{\bar{S_1}}((1,0,0))\geq |F^*||K^*|$).   So we have the following cases:\\
\textbf{Case 1} $|F|\geq 4$ and $|K|\geq 4$: In this cases we get a contradiction.\\
\textbf{Case 2} $|F|=|K|=3$: The zero-divisor graph of $\Z_2 \times F \times K$ is  illustrated  in Figure \ref{fig4}.
 \begin{figure}[ht]
	\centering
	\includegraphics[scale=0.5]{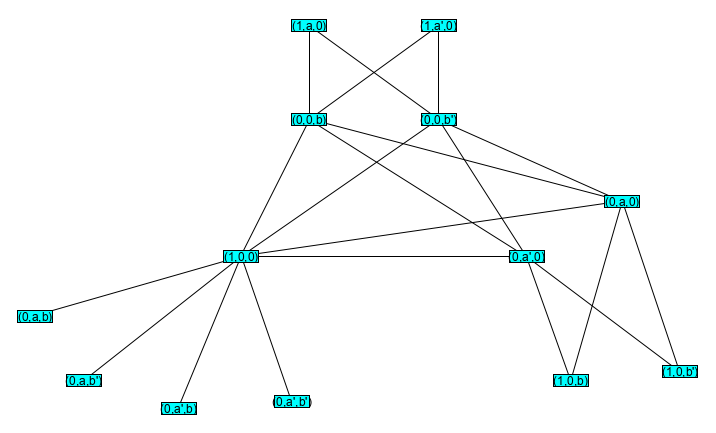}
	\caption{$\Gamma(\Z_2\times F\times K)$}
	\label{fig4}
\end{figure}   
Then, $S_1$ contains at least three vertices from the set $\{(0,0,b), (0,0,b'), (0,a,0),(0,a',0)\}$, so we can assume that the $\{(0,0,b), (0,0,b'), (0,a,0)\}\subset  S_1$. Then, $\{(1,a,0),(1,a',0)\}\subset S_2$ (since $S_2$ is a dominating set). Thus, $deg_{S_2}((1,a,0))+1\geq deg_{\bar{S_2}}((1,a,0))$ and so $1\geq 2$, a contradiction.\\
\textbf{Case 2} $|F|=4$ and $|K|=3$:  The zero-divisor graph of $\Z_2 \times F \times K$ is  illustrated  in Figure \ref{fig5}.  
 \begin{figure}[ht]
	\centering
	\includegraphics[scale=0.5]{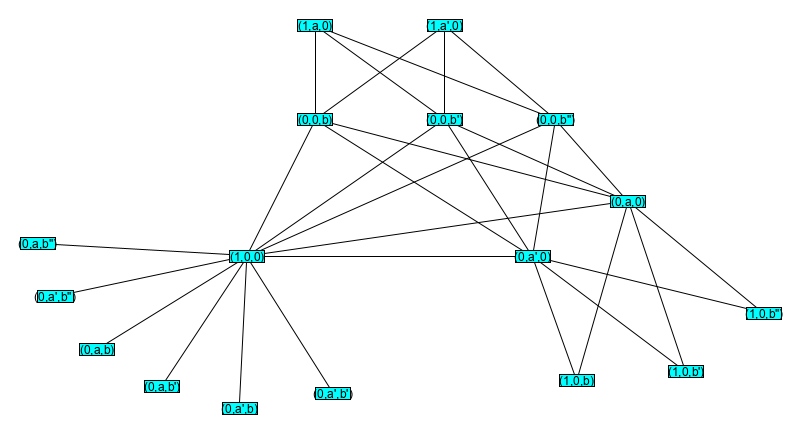}
	\caption{$\Gamma(\Z_2\times F\times K)$}
	\label{fig5}
\end{figure} 
Then,  $\{(0,0,b), (0,0,b'), (0,0,b''), (0,a,0),(0,a',0)\}\subset S_1 $ (since $deg_{S_1}((1,0,0))+1\geq deg_{\bar{S_1}}((1,0,0))$ ). Thus,  $deg_{S_2}((1,a,0))+1\geq deg_{\bar{S_2}}((1,a,0))$ implies that $1\geq 3$, a contradiction.\\
Hence, $\psi_g(\Gamma(\Z_2\times F \times K))<2$ and so $\Gamma(\Z_2\times F \times K)$ is not  partitionable into global defensive alliance.
\cqfd

\begin{cor}
	Let $p,q\geq 3$ be two  prime numbers. Then, $\psi_g(\Gamma(\Z_2\times \Z_p\times \Z_q))=1$.
\end{cor}

We end this paper by investigating rings with small global defensive alliance number and small global defensive alliance partition number. Namely,  $\gamma_{a}(\Gamma(R))=1,2$ and $\psi_{g}(\Gamma(R))=2,3$.

\begin{thm}\label{thm_partitioning1}
	Let $R$ be a finite ring such that $\gamma_a(\Gamma(R))=1$. Then, the following statements are equivalent:
	\begin{enumerate}
		\item  $\Gamma(R)$ is partitionable into global defensive alliances.
		\item  $|Z(R)|=3$.
		\item  $R$ is isomorphic to one of the rings $\Z_2\times \Z_2$, $\Z_9$, $\Z_3[X]/(X^2)$.
	\end{enumerate}	 
\end{thm}
\pr $(1)\Rightarrow (2)$ Assume that  $\Gamma(R)$ is partitionable into global defensive alliances, then $\psi_g(\Gamma(R))\geq 2$ and so $|Z(R)|=3$, by Proposition \ref{prop_|Z(R)|} and \cite[Proposition 2.2]{NA}.\\
$(2)\Rightarrow (3)$ Follows from \cite[Corollary 1]{MBRB}.\\
$(3)\Rightarrow (1)$ The zero-divisor graphs  of these rings are isomorphic to a simple graph with two vertices  and one edge. Hence, we get the result.
\cqfd

\begin{cor}
	Let $R$ be a finite ring such that $\Gamma(R)$ is partitionable into global defensive alliances. If $\gamma_a(\Gamma(R))=1$,  then $\psi_g(\Gamma(R))=2$.
\end{cor}

\begin{thm}\label{thm_partitioning2}
	Let $R$ be a finite ring such that $\gamma_a(\Gamma(R))=2$. Then, $\Gamma(R)$ is partitionable into global defensive alliances if and only if either $\psi_g(\Gamma(R))=2$ or  $R\cong \F_4\times \F_4$.
\end{thm}
\pr $\Rightarrow)$  Assume that  $\Gamma(R)$ is partitionable into global defensive alliances and suppose that $R\ncong \F_4\times \F_4$. Then $\psi_g(\Gamma(R))\geq 2$.  So, we need just to prove the other equality (i.e., $\psi_g(\Gamma(R))\leq 2$ ).  We have $\gamma_a(\Gamma(R))=2$, then by Proposition \ref{prop_|Z(R)|} and \cite[Proposition 2.2]{NA}, $\psi_g(\Gamma(R))^2-\psi_g(\Gamma(R))-6\leq 0$ and so  $\psi_g(\Gamma(R))$ is either $2$ or $3$. Suppose that $\psi_g(\Gamma(R))=3$, then by Proposition \ref{prop_|Z(R)|}, $|Z(R)|\geq 7$. Thus,   using  \cite[Proposition 3.3]{NA}, $R\cong \F_4\times \F_4$, a contradiction by the hypothesis. Thus,  $\psi_g(\Gamma(R))=2$.\\
$\Leftarrow)$ Its obvious.
\cqfd

\begin{cor}\label{cor_partition2}
	Let $R$ be a finite ring such that $\gamma_a(\Gamma(R))=2$. Then, we have the following equivalents:
	\begin{enumerate}
	\item $\psi_g(\Gamma(R))=2 $ if and only if $R$ is isomorphic to one of the rings $\Z_2\times \Z_4$, $\Z_2\times \Z_{2}[X]/(X^2)$, $\Z_3\times \Z_3$, $\Z_3\times \F_4$, $\Z_{25}$ and $\Z_5[X]/(X^2)$. 
	\item $\psi_g(\Gamma(R))=3$ if and only if $R\cong \F_4\times \F_4$.
	\end{enumerate}    
\end{cor}
\pr It follows from Theorem \ref{thm_partitioning2} and  \cite[Proposition 3.3]{NA}. \cqfd


\bigskip
		
		
		
\begin{thebibliography}{999}\addcontentsline{toc}{section}{\protect\numberline{}{Bibliography}}
	



\bibitem{SAAM} S. Akbari and A. Mohammadian,   On the zero-divisor graph of a commutative ring, \textit{ J. Algebra } \textbf{274} (2004)  847–855.

			
	
\bibitem{SHS}	S. Akbari, H. R. Maimani and S. Yassemi, When a zero-divisor graph is planar or a complete $r$-partite graph, \textit{J. Algebra} \textbf{270} (2003) 169–180.
	
	
\bibitem{AFLL01}  D. F. Anderson,   A.  Frazier, A. Lauve and  P. S.  Livingston, The Zero-Divisor Graph of a Commutative ring, II, \textit{ Lect. Notes Pure Appl. Math.}   \textbf{220} (2001)  61-72.
	
	
\bibitem{ADL} D. F. Anderson and  P. S. Livingston,  The zero-divisor graph of a commutative ring, \textit{ J. Algebra} \textbf{217} (1999) 434-447.
	
	
	
\bibitem{DM} D. D. Anderson and   M. Naser,  Beck's coloring of commutative ring, \textit{J. Algebra} \textbf{159} (1993)  500-514.
	
	
	
	
\bibitem{A08} D. F. Anderson,   On the diameter and girth of a zero-divisor graph II, \textit{ Houston J. Math}  \textbf{34}    (2008) 361-371.  
	
	
	
\bibitem{BAD} D. F. Anderson and A. Badawi,   On the zero-divisor graph of a ring, \textit{ Comm. Algebra  } \textbf{36} (2008) 3073-3092.
	
	
	
	
	
	
	
\bibitem{AM07} D. F. Anderson and S. B.  Mulay,  On the diameter and girth of a zero-divisor graph, \textit{ J. Pure Appl. Algebra} \textbf{210} (2007) 543-550.


	
	
	
\bibitem{COY} M. Axtell, J.  Coykendall and  J.  Stickles,  Zero-divisor graphs of polynomial and power series over commutative rings, \textit{ Comm. Algebra } \textbf{33} (2005) 2043-2050.
	
	
	
	
	
	
	
	
\bibitem{BI} I. Beck,  Coloring of commutative rings, \textit{ J. Algebra} \textbf{116} (1988) 208-226.
	
	
	
	
	
\bibitem{MBRB} M. Behboodi and R. Beyranvand,  On the structure of commutative rings with $p^{k_1}_1\ldots p^{k_n}_n (1 \leq k_i \leq 7)$  zero-divisors,  \textit{ Eur. J. Pure Appl. Math.} \textbf{3} (2010) 303–316.		
			




\bibitem{DBK} D. Bennis, B. E. Alaoui and K. Ouarghi,   On global  defensive k-alliances in zero-divisor graphs of finite commutative rings, \textit{ J. Algebra  Appl.} (2022) 2350127.	






\bibitem{RGOJ} R. J. Morales, G. R. Hernandez, O. R. Cayetano and J. R. Valencia,  On global offensive alliance in zero-divisor graphs, \textit{Mathematics} \textbf{10} (2022) 298.	
	
	
	
	
\bibitem{HHH}  T. W. Haynes, S. T. Hedetniemi, M. A. Henning, \textit{ Topics in domination in graphs},  (Developments in Mathematics, Springer, Cham. \textbf{64} 2020).
	
	
	
	
\bibitem{TSP} T. W. Haynes, S. T.  Hedetniemi and P. J. Slater,  \textit{ Fundamentals of Domination in Graphs} ( Marcel Dekker, Inc. New York,  1998).
	
	
	
	
\bibitem{TSP2} T. W. Haynes, S. T. Hedetniemi and P. J. Slater,  \textit{ Domination in Graphs: Advanced Topics} ( Marcel Dekker, Inc. New York, 1998).
	
	
	
\bibitem{TSM} T. W. Haynes, S. T. Hedetniemi and M. A.  Henning,  Global defensive alliances in graphs,  \textit{Electron. J. Combin.} \textbf{10} (2003) R47.	
	
	
	
	
\bibitem{PSS} P. Kristiansen, S. M. Hedetniemi and  S. T. Hedetniemi,  Alliances in Graph, \textit{J. Combin. Math. Combin. Comput.} \textbf{177}  (2004) 48-155.
	
	
	
	
	
\bibitem{PKSMST} P. Kristiansen, S.M. Hedetniemi and S.T. Hedetniemi, Introduction to alliances in graphs, \textit{IBID} (2002) 308-312.	
	






\bibitem{NA} N. Muthana and A. Mamouni,   On defensive alliance in zero-divisor graphs, \textit{ J. Algebra  Appl.} \textbf{20}  (2021)  2150155.	

	
	
	
	
\bibitem{JJ} J.A.  Rodr\'iguez-Vel\'azquez and  J.M. Sigarreta,  Global defensive $k$-alliances in graphs, \textit{ Discrete Appl. Math.} \textbf{157} (2009) 211–218.
	
	
	
	
\bibitem{JIJ} J.A. Rodr\'iguez-Vel\'azquez, I.G. Yero, and J.M.  Sigarreta,   Defensive $k$-alliances in graphs, \textit{ Appl. Math. Lett.}   \textbf{22}  (2009)  96-100.  
	
	
	
	
\bibitem{Shaf} K. H. Shafique,  \textit{ Partitioning a graph in alliances and its application to data clustering} ( PhD thesis, School of Computer Science, University of Central Florida, Orlando, FL, 2001).
	




\bibitem{KR} K.H. Shafique and R.D.  Dutton,  On satisfactory partitioning of graphs, \textit{ Congr. Numer.}  \textbf{154} (2002)  183–194.	





		
		
\bibitem{NOS}	N. O. Smith, Planar zero-divisor graphs, Focus on Commutative Rings Research (Nova Science Publisher, New York, 2006), pp. 177–186.		

	
	
		
		
		
			
\bibitem{ISJJ} I.G. Yero, S. Bermudo, J.A.  Rodríguez-Velázquez and  J.M.  Sigarreta,   Partitioning a graph into defensive $k$-alliances,  \textit{Acta Math. Sin. (Engl. Ser.) } \textbf{27} (2011) 73–82.	
		
		
		
	
		
			
			
\end{thebibliography}
	\end{document}